\documentclass[11pt]{amsart}

\usepackage{amssymb}
\usepackage{amsmath}
\usepackage{enumerate}
\usepackage[latin1]{inputenc}
\usepackage[dvips]{graphicx}

\def\s{\mathbb{S}}
\def\h{\mathbb{H}}
\def\r{\mathbb{R}}

\def\c{\mathbb{C}}
\def\p{\mathbb{P}}

\newtheorem{remark}{Remark}
\newtheorem{theorem}{Theorem}
\newtheorem{proposition}{Proposition}
\newtheorem{corollary}{Corollary}

\numberwithin{equation}{section}

\DeclareMathOperator{\arctanh}{arctanh}
\DeclareMathOperator{\arccosh}{arccosh}

\begin{document}

\title[Hamiltonian stationary Lagrangian spheres]{On Hamiltonian stationary Lagrangian spheres
in non-Einstein K\"{a}hler surfaces}

\author{Ildefonso Castro}
\address{Departamento de Matem\'{a}ticas \\
Universidad de Ja\'{e}n \\
23071 Ja\'{e}n, SPAIN} \email{icastro@ujaen.es}

\author{Francisco Torralbo}
\address{Departamento de Geometr\'{\i}a  y Topolog\'{\i}a \\
Universidad de Granada \\
18071 Granada, SPAIN} \email{ftorralbo@ugr.es}

\author{Francisco Urbano}
\address{Departamento de Geometr\'{\i}a  y Topolog\'{\i}a \\
Universidad de Granada \\
18071 Granada, SPAIN} \email{furbano@ugr.es}

\thanks{Research partially supported by a MEC-Feder grant MTM2007-61775}



\date{}

\begin{abstract}
Hamiltonian stationary Lagrangian spheres in K\"ahler-Eins- tein
surfaces are minimal. We prove that in the family of non-Einstein
K\"ahler surfaces given by the product $\Sigma_1\times\Sigma_2$ of
two complete orientable Riemannian surfaces of different constant
Gauss curvatures, there is only a (non minimal) Hamiltonian
stationary Lagrangian sphere. This example is defined when the
surfaces $\Sigma_1$ and $ \Sigma_2$ are spheres.
\end{abstract}

\maketitle

\section{Introduction}

Minimal Lagrangian submanifolds of K\"ahler manifolds acquired a
great importance since 1982, when Harvey and Lawson in \cite{HL}
proved that the orientable minimal Lagrangian submanifolds of a
Calabi-Yau manifold are exactly the special Lagrangian
submanifolds, i.e.\ the calibrated submanifolds for the special
Lagrangian calibration that a Calabi-Yau manifold carries. In
particular, they are volume minimizing in their integer homology
classes.

For Lagrangian submanifolds, it is also natural to study those
which are critical for the volume under variations that keep the
property to be Lagrangian. This kind of submanifolds are called
{\it Lagrangian stationary} and  Schoen and Wolfson \cite{SW}
proved that, when the ambient K\"ahler manifold is Einstein, they
are exactly the minimal Lagrangian ones. Among these variations,
the {\it Hamiltonian variations} are particulary important and Oh
was who studied in \cite{O} Lagrangian submanifolds of a K\"ahler
manifold which are critical for the volume under Hamiltonian
deformations. These submanifolds are called {\it Hamiltonian
stationary} (or also {\it Hamiltonian minimal}) and the
corresponding Euler-Lagrange equation that characterizes them is
\begin{equation}\label{EL}
d^*\alpha_{H}=0,
\end{equation}
where $d^*$ is the adjoint operator of the exterior differential
and $\alpha_H$ is the $1$-form on the submanifold defined by
$\alpha_{H}=H\lrcorner\,\omega$, $\omega$ being the K\"ahler
$2$-form of the ambient space and $H$ the mean curvature vector of
the submanifold. It is clear that this condition is equivalent to
$\hbox{div}\,JH=0$, where $J$ is the complex structure of the
ambient space and div stands for the divergence operator on the
submanifold.

It is a very well known fact that if $\Sigma$ is a Lagrangian
submanifold of a K\"ahler manifold, then
\begin{equation}\label{Ric}
d\alpha_H=\hbox{Ric}\,_{|\Sigma},
\end{equation}
where Ric denotes the Ricci $2$-form on the ambient space. Hence
if the ambient space is K\"ahler-Einstein, then $\alpha_H$ is a
closed $1$-form. In this case, using (\ref{EL}) it is clear that
if $\Sigma $ is Hamiltonian stationary then $\alpha_H$ is a
harmonic $1$-form on $\Sigma$. This fact implies that if $\Sigma$
is, in addition,  compact and its first Betti number vanishes,
then $H=0$, and the submanifold is minimal. Thus, in particular,
{\it Hamiltonian stationary Lagrangian spheres in
K\"ahler-Einstein manifolds are necessarily minimal}.

Most of the literature about Hamiltonian stationary Lagrangian
surfaces is centered in complex space forms. Hence when the
surface is compact, the study of the Hamiltonian stationary
Lagrangian tori plays a relevant role (for example, see
\cite{CU1}, \cite{HR1} or \cite{HR2}).

In this paper we are interested in studying non minimal
Hamiltonian stationary spheres of K\"ahler manifolds. As the
ambient space can not be Einstein, for the sake of simplicity we
have considered (as K\"ahler target manifolds) the products of two
Riemann surfaces $\Sigma_1\times\Sigma_2$ endowed with metrics of
constant Gauss curvatures $c_1$ and $c_2$ with $c_1\not= c_2$,
because $\Sigma_1\times\Sigma_2$ is Einstein if and only if
$c_1=c_2$. It is clear that we can assume without restriction that
$c_1>c_2$.

Our main contribution in this paper is the following Existence and
Uniqueness Theorem.

\begin{theorem}\label{Th1}
Let $\Phi:\Sigma \rightarrow \Sigma_1 \times \Sigma_2$ be a
Hamiltonian stationary Lagrangian immersion of a sphere $\Sigma$
in the product of two complete orientable Riemannian surfaces
$\Sigma_i$, $i=1,2$, of constant Gauss curvatures $c_i$, $i=1,2$,
with $c_1 > c_2$. Then
\begin{enumerate}
\item $c_1, c_2 >0 $ and hence $\Sigma_1\times\Sigma_2$ is isometric to the
product of two spheres $\s^2_{c_1}\times\s^2_{c_2}$.
\item $\Phi$ is
congruent to the Hamiltonian stationary Lagrangian immersion
$\Phi_0 =(\phi_0,\psi_0): \s^2 \rightarrow \s^2_{c_1} \times
\s^2_{c_2}$, where
\begin{equation}\label{Phi_0}
\begin{array}{c}
 \phi_0(z,x)=\frac{\textstyle
2\sqrt{c_1-c_2}}{\textstyle c_1 - c_2 x^2} \left( ixz,
\frac{\textstyle -c_1+(2c_1-c_2)x^2}{\textstyle 2\sqrt{c_1}\sqrt{c_1-c_2}} \right), \\ \\
\psi_0(z,x)=\frac{\textstyle 2\sqrt{c_1-c_2}}{\textstyle c_1 - c_2
x^2} \left(
 \bar{z},
\frac{\textstyle (c_1-2c_2)+c_2 x^2}{\textstyle
2\sqrt{c_2}\sqrt{c_1-c_2}} \right).
\end{array}
\end{equation}
\end{enumerate}
$\s^2$ is the unit sphere and, in general, the sphere $\s^2_c$ is
given by $\s^2_c=\{(z,x)\in\c\times\r\,/\,|z|^2+x^2=1/c\}$.
\end{theorem}

\begin{remark}
{\rm It is an easy exercise to check that $\Phi_0$ is a
Hamiltonian stationary Lagrangian embedding except at the poles of
$\s^2$ where it has a double point. Hence
$\mathcal{S}_0=\Phi_0(\s^2)$, which is given by
\[
\begin{array}{c}
\mathcal{S}_0 = \{ ((z_1,x_1),(z_2,x_2))\in \s^2_{c_1} \!\times\!
\s^2_{c_2}\,/\, \Re (z_1 z_2)= 0, \,  \sqrt{c_2} x_1
\!-\!\sqrt{c_1} x_2 =\frac{ c_2\!-\!c_1}{\sqrt{c_1c_2}} \},
\end{array}
\]
is a surface of $\s^2_{c_1} \!\times\!\s^2_{c_2}$ with a
singularity at $((0,1/\sqrt{c_1}),(0,1/\sqrt{c_2}))$.

In the hypothesis of Theorem \ref{Th1}, we have not assumed that
the submanifold was non minimal, because in these ambient spaces
there are no minimal Lagrangian spheres (see Proposition
\ref{Prop2}).}
\end{remark}

The proof of Theorem \ref{Th1} involves arguments of different
nature which are developed in section 3.  In section 3.1 we show
that the property of the surface to be Hamiltonian stationary,
joint with the product structure of the ambient space, allow to
define a holomorphic Hopf differential (see (\ref{Theta})) on the
surface. A similar construction was made by Abresh and Rosenberg
\cite{AR} for constant mean curvature surfaces in $\s^2\times\r$
and $\h^2 \times \r$, where $\h^2$ is the hyperbolic plane. As the
surface is a sphere, this Hopf differential vanishes, and this
allows us to construct another holomorphic Hopf differential (see
(\ref{Xi})), which can be defined only for Hamiltonian stationary
Lagrangian spheres.

The vanishing of the two Hopf differentials provide, in section
3.2, important geometric properties of our surface. First we prove
that the length of the mean curvature $|H|$ is locally the product
of a smooth positive function and the modulus of a holomorphic
function. Hence we have that
\[
\int_{\Sigma}\Delta\,\hbox{log}\,|H|\,dA=-2\pi N(|H|),
\]
where $N(|H|)$ is the sum of the orders of the zeroes of $|H|$ and
$\Delta$ is the Laplacian operator. Second we show how this
formula becomes in the fundamental equality
\begin{equation}\label{fla}
\int_{\Sigma}\left(|H|^2+\frac{c_1+c_2}{4}\right)dA=8\pi.
\end{equation}
This fact, joint with the Bochner formula for the tangent vector
field $JH$, allows to prove that $c_1$ and $c_2$ must be positive
and so $\Sigma_1\times\Sigma_2=\s^2_{c_1}\times\s^2_{c_2}$.

Finally, in section 3.3, if we consider
$\s^2_{c_1}\times\s^2_{c_2}\subset\r^6$ and $\tilde{H}$ is the
mean curvature of $\Sigma$ into $\r^6$, then the equality
(\ref{fla}) give us the value of its Willmore functional:
\[
W(\Phi):=\int_{\Sigma}|\tilde{H}|^2dA=8\pi.
\]
As for topological reasons our surface can not be embedded (see
Corollary 1), a result of Li and Yau \cite{LY} yields that our
surface is not only a Willmore surface in $\r^6$ but also a
minimizer for the Willmore functional. Moreover, our surface is
the compactification by an inversion of $\r^6$ of a complete
minimal surface of $\r^6$ with two planar embedded ends and total
curvature $-4\pi$, which turns out to lie in a $4$-dimensional
affine subspace of $\r^6$. Using the classification obtained by
Hoffman and Osserman in \cite{HO}, we obtain a $1$-parameter
family of Lagrangian immersions $\Phi_{t},\, t\geq 0$ (see
(\ref{Phis}) for their explicit definition) of a sphere into
$\s^2_{c_1}\times\s^2_{c_2}$ which satisfy (\ref{fla}). Next it is
easy to check that the immersion $\Phi_0$ is the only Hamiltonian
stationary in this family. We remark that precisely $\Phi_0$ is
the compactification of the Lagrangian catenoid defined in
\cite{HL}.

The above reasoning also proves the following result.

\begin{theorem}\label{Th2}
Let $\Phi:\Sigma \rightarrow \s^2_{c_1} \times \s^2_{c_2}$ be a
Lagrangian immersion of a sphere $\Sigma$ in the product of two
spheres of Gauss curvatures $c_1 > c_2$. Then
\[
\int_\Sigma \left( |H|^2 + \frac{c_1+c_2}{4} \right)dA \geq 8\pi ,
\]
and the equality holds if and only if $\Phi$ is congruent to some
$\Phi_{t},\, t\geq 0$ (see (\ref{Phis}) for their definition).
\end{theorem}

The Lagrangian variation $\Phi_{t}$ of $\Phi_0$ is Hamiltonian,
and if $A(t)$ is the area of the induced metric on $\Sigma$ by
$\Phi_{t}$, then $A(t)$ has a global maximum at $t=0$ and
$\lim_{t\rightarrow\infty}A(t)=0$. Therefore {\it $\Phi_0$ is
unstable and there are no minimizers for the area in its
Hamiltonian isotopy class}.

\section{Preliminaries}

Let $(\Sigma_1,g_1)$ and $(\Sigma_2,g_2)$ be two Riemannian
orientable surfaces with constant Gauss curvatures $c_1$ and
$c_2$. We denote by $\omega_1$ and $\omega_2$ their K\"{a}hler
$2$-forms and by $J^1$ and $J^2$ their complex structures as
Riemann surfaces.

Consider the product K\"ahler surface  $\Sigma_1\times \Sigma_2$,
endowed with  the product metric $\langle,\rangle=g_1\times g_2$
and the complex structure $J=(J^1,J^2)$. Its K\"{a}hler $2$-form
is $\omega=\pi_1^*\omega_1+\pi_2^*\omega_2$, where $\pi_i$ are the
projections of $\Sigma_1\times \Sigma_2$ onto $\Sigma_i$, $i=1,2$.

Along the paper we will use the {\em product structure} on
$\Sigma_1\times\Sigma_2$, which is the (1,1)-tensor $P$ defined by
\[ P:T(\Sigma_1\times \Sigma_2)\rightarrow
T(\Sigma_1\times \Sigma_2), \  P(v_1,v_2)=(-v_1,v_2).
\]
We remark that $P$ is parallel, i.e.\ $\overline{\nabla}P=0$,
where $\overline{\nabla}$ is the Levi-Civita connection of the
product metric $\langle, \rangle$. In addition,  $P$ is an
$\langle , \rangle$-isometry and verifies $P^2=$Id and $P\circ J
=J\circ P$.

\vspace{0.3cm}

 Let
$\Phi=(\phi,\psi):\Sigma\rightarrow\Sigma_1\times \Sigma_2$ be a
Lagrangian immersion of a surface $\Sigma$, i.e.\ an immersion
satisfying $\Phi^*\omega=0$. In this case, this condition means
that $\phi^*\omega_1+\psi^*\omega_2=0$ or, equivalently:
\[
0=\langle Jd\Phi_p(v),d\Phi_p(w)\rangle = g_1( J^1
d\phi_p(v),d\phi_p(w))+ g_2( J^2 d\psi_p(v),d\psi_p(w)),
\]
for any $p\in\Sigma $ and $v,w\in T_p\Sigma$. If $\Sigma$ is an
oriented surface and $\omega_{\Sigma}$ is the area $2$-form of its
induced metric, the Jacobians of $\phi$ and $\psi$ are given by
\begin{equation}\label{Jacs}
\phi^*\omega_1=\hbox{Jac}\,(\phi)\,\omega_{\Sigma},\quad\psi^*\omega_2=\hbox{Jac}\,(\psi)\,\omega_{\Sigma},
\end{equation}
and hence $\hbox{Jac}\,(\phi)=-\hbox{Jac}\,(\psi)$. We will call
the function
\begin{equation}\label{defC}
C:=\hbox{Jac}\,(\phi)=-\hbox{Jac}\,(\psi)
\end{equation}
the {\em associated Jacobian} of the oriented Lagrangian surface
$\Sigma$. It is easy to check that the restriction to $\Sigma$ of
the Ricci $2$-form on $\Sigma_1\times\Sigma_2$ is given by
$\hbox{Ric}\,_{|\Sigma}=\frac{1}{2}(c_1-c_2)C\omega_{|\Sigma}$ and
then (\ref{Ric}) becomes in
\begin{equation}\label{d_Maslov}
d\alpha_{H}=\frac{1}{2}(c_1-c_2)C\omega_{|\Sigma}.
\end{equation}
In the compact case, integrating (\ref{d_Maslov}) and using the
Stokes theorem and (\ref{Jacs}) and (\ref{defC}), we get the
following result. \vspace{0.1cm}

\begin{proposition}\label{Prop1}
If $\Phi=(\phi,\psi):\Sigma \rightarrow \Sigma_1 \times \Sigma_2$
is a Lagrangian immersion of a compact oriented surface $\Sigma$
in the product $\Sigma_1\times\Sigma_2$ with $\Sigma_1$ and
$\Sigma_2$ compact and $c_1\not= c_2$, then the degrees of $\phi$
and $\psi$ satisfy
\[
\mathrm{deg}\,(\phi)=\mathrm{deg}\,(\psi)=0.
\]
\end{proposition}
If $\Sigma_1\times\Sigma_2=\s^2_{c_1} \times \s^2_{c_2}$ this
proposition offers an interesting consequence.
\begin{corollary}
If $\Phi=(\phi,\psi):\Sigma \rightarrow \s^2_{c_1} \times
\s^2_{c_2}$ is a Lagrangian embedding of a compact orientable
surface $\Sigma$, with $c_1\neq c_2$, then $\Sigma$ must be a
torus.
\end{corollary}
\begin{remark}
{\rm  This result is not true when $c_1=c_2$ since the graph of the antipodal map is a Lagrangian sphere embedded
in $\s^2_{c_1} \times
\s^2_{c_1}$ (see \cite{CU2}).}
\end{remark}
\begin{proof}
Let $\gamma $ be the genus of $\Sigma $. Following the same
reasoning that in the proof of Proposition 3 in \cite{CU2} (that
is, computing the self-intersection number of $\Phi$), as $\Phi$
is an embedding, we can prove that the degrees of $\phi$ and
$\psi$ satisfy
$$\gamma=1+\mathrm{deg}\phi\,\mathrm{deg}\psi.$$
Proposition
\ref{Prop1} says that $\gamma =1$ and so $\Sigma $ is a torus.
\end{proof}

In the general case ($\Sigma$ not necessarily compact), we can
prove the following result. \vspace{0.1cm}

\begin{proposition}\label{Prop2}
Let $\Phi:\Sigma \rightarrow \Sigma_1 \times \Sigma_2$ be a
Lagrangian immersion of an orientable surface $\Sigma$ in the
product $\Sigma_1\times\Sigma_2$ with $c_1\not= c_2$. If $\Phi $
is minimal (or, more generally, $\Phi$ has parallel mean curvature
vector), then $\Phi $ must be a product of two curves in
$\Sigma_1$ and $\Sigma_2$ with constant curvatures. In particular,
there are no minimal Lagrangian spheres in
$\Sigma_1\times\Sigma_2$ with $c_1\not= c_2$.
\end{proposition}
\begin{proof} As usual, $\Phi=(\phi,\psi)$. Using again (\ref{d_Maslov}),
the assumption about the mean curvature joint with $c_1\not= c_2$
say that the associated Jacobian $C\equiv 0$. But it is not
difficult to check (see Lemma 2.1 in \cite{CU2}) that if
$\{e_1,e_2\}$ is an oriented orthonormal basis of $T_p\Sigma$,
then
\begin{equation}\label{1=bo}
|d\phi_p(e_1)|^2+|d\phi_p(e_2)|^2=|d\psi_p(e_1)|^2+|d\psi_p(e_2)|^2=1,
\end{equation}
which implies that the ranks of $d\phi$ and $d\psi$ at any point
of $\Sigma $ are always positive. Hence, if $C\equiv 0$, i.e.\ the
Jacobians of $\phi$ and $\psi$ vanish, then the ranks of $d\phi$
and $d\psi$ at any point must be necessarily $1$ and so both
functions $\phi$ and $\psi$ define curves in $\Sigma_i$, $i=1,2$,
respectively. Now it is an exercise to check that the product of
two curves has parallel mean curvature vector if and only if both
curves have constant curvature.
\end{proof}

Finally, if  $\bar{R}$ denotes the curvature operator of
$\Sigma_1\times \Sigma_2$, it is easy to prove that
$\bar{R}(e_1,e_2,e_2,e_1)=(c_1+c_2)C^2$, where $\{e_1,e_2\}$ is an
orthonormal frame on $T\Sigma$. Thus the Gauss equation of $\Phi$
can be written as
\begin{equation}\label{EqG}
K=(c_1+c_2)C^2+2|H|^2-\frac{|\sigma|^2}{2},
\end{equation}
where $K$ is the Gauss curvature of $\Sigma$, $H$ the mean
curvature of $\Phi$ and $\sigma$ the second fundamental form of
$\Phi$.

\section{Proof of the Main Theorem}

Let $\Phi:\Sigma \rightarrow \Sigma_1 \times \Sigma_2$ be a
Hamiltonian stationary Lagrangian immersion of a sphere $\Sigma$
in the product of two orientable Riemannian complete surfaces
$\Sigma_i$ of constant Gauss curvatures $c_i$, $i=1,2$, with $c_1
> c_2$.

We split the proof of the Theorem \ref{Th1} in three parts.

\subsection{Hopf differentials}

We consider a local complex coordinate $z=x+iy$ on $\Sigma$ and
let $\partial_z=\frac{1}{2}\left (\frac{\partial}{\partial
x}-i\frac{\partial}{\partial y}\right )$ and $\partial_{\bar
z}=\frac{1}{2}\left (\frac{\partial}{\partial
x}+i\frac{\partial}{\partial y}\right )$ the Cauchy-Riemann
operators. Letting $g$ and $\nabla $ also denote the complex
extension of the induced metric and its Levi-Civita connection, we
have that
\[
\nabla_{\partial_z}\partial_z =2u_z \partial_z, \quad
\nabla_{\partial_z}\partial_{\bar z} =0,
\]
where $g=e^{2u}|dz^2|$.

We will also use $\nabla $ to denote the Levi-Civita connection on
the pulled-back bundle $\Phi^{-1} T(\Sigma_1\times \Sigma_2)$ and
let
\[
\delta = \nabla_{\partial_z}, \quad \bar{\delta} =
\nabla_{\partial_{\bar z}}
\]
be the corrresponding covariant derivatives acting on $\Phi^{-1}
T(\Sigma_1\times \Sigma_2)$.

In the sequel we will frequently identify the fibres of $\Phi^{-1}
T(\Sigma_1\times \Sigma_2)$ and $T(\Sigma_1\times \Sigma_2)$ in
order to write
\[
\delta \Phi= \Phi_* (\partial_z), \quad \bar{\delta} \Phi= \Phi_*
(\partial_{\bar z}).
\]
Hence, we have that
\begin{equation}\label{conformal}
\langle \delta \Phi,\delta \Phi \rangle = 0,\quad |\delta
\Phi|^2=e^{2u}/2.
\end{equation}
From (\ref{defC}) it is easy to get that the associated Jacobian
$C$ can be written as
\begin{equation}\label{C}
C =
-i\,e^{-2u}\langle P\delta \Phi,
J\bar{\delta}\Phi \rangle.
\end{equation}
Along the proof it will be useful the following formula deduced
from (\ref{1=bo}), (\ref{conformal}) and (\ref{C}):
\begin{equation}\label{Frenet3}
P\delta \Phi = 2 (e^{-2u}\langle \delta \Phi, P\delta \Phi \rangle
\bar{\delta}\Phi + i\,C\, J\delta \Phi).
\end{equation}

The starting point in the proof is to use (\ref{d_Maslov}) and
(\ref{C}) to get the fundamental relation:
\[
 \partial_{\bar z}\langle H, J\delta \Phi \rangle = -\frac{e^{2u}}{4}\left(
{\mathrm{div} JH} - \frac{i}{2} (c_1-c_2)C \right).
\]
Hence {\it $\Phi $ is Hamiltonian stationary if and only if
\begin{equation}\label{h_z}
\partial_{\bar z}\langle H, J\delta \Phi \rangle = \frac{i}{8}
(c_1-c_2)e^{2u}C.
\end{equation}}

We use (\ref{h_z}) and the product structure on
$\Sigma_1\times\Sigma_2$ to define a Hopf quadratic differential
on $\Sigma$. Since $H=2e^{-2u}\bar{\delta}\delta\Phi$, using the
properties of the product structure $P$, we deduce that
\[
\partial_{\bar z}\langle \delta \Phi,P \delta\Phi \rangle =
\langle\bar{\delta}\delta\Phi, P\delta \Phi
\rangle+\langle\delta\Phi,P\bar{\delta}\delta\Phi\rangle =e^{2u}
\langle H,P\delta\Phi\rangle.
\]
We observe that if $\xi$ is a normal vector field, it is easy to
check that $\langle \xi, P \delta \Phi \rangle =2i\,C \langle \xi,
J \delta \Phi \rangle$ using (\ref{C}) and (\ref{Frenet3}). So
$\langle H, P \delta \Phi \rangle =2i\,C \langle H, J \delta \Phi
\rangle$. Putting this in the above equation leads to
\begin{equation}\label{G_zbarra}
\partial_{\bar z}\langle \delta \Phi,P \delta\Phi \rangle
=2i\,e^{2u}C\langle H, J\delta \Phi \rangle .
\end{equation}
Now (\ref{h_z}) and (\ref{G_zbarra}) yield to $
\partial_{\bar z}\left( \langle \delta \Phi,P \delta\Phi \rangle -
\frac{8}{c_1-c_2} \langle H, J\delta \Phi \rangle ^2 \right)=0 $.
We have proved that {\it the Hopf quadratic differential
\begin{equation}\label{Theta}
\Theta(z)=\left( \langle \delta \Phi,P \delta\Phi \rangle -
\frac{8}{c_1-c_2} \langle H, J\delta \Phi \rangle ^2 \right)
(dz)^2
\end{equation}
is holomorphic.}

As the surface $\Sigma $ is a sphere, $\Theta\equiv 0$ and so we
get that
\begin{equation}\label{ya_esfera}
 \langle \delta \Phi,P\delta \Phi \rangle =
\frac{8}{c_1-c_2} \langle H, J\delta \Phi \rangle ^2 .
\end{equation}
It is easy to check that
\begin{equation}\label{|G|y|h|}
 |\langle \delta \Phi, P\delta \Phi
\rangle|^2 =\frac{e^{4u}}{4}(1-4C^2),\quad  |\langle H,
J\delta  \Phi \rangle|^2=\frac{e^{2u}}{4}|H|^2 ,
\end{equation}
and, as a consequence, taking modules in (\ref{ya_esfera}) and using
(\ref{|G|y|h|}), we obtain
\begin{equation}\label{|H|oC}
|H|^2=\frac{c_1-c_2}{4}\sqrt{1-4C^2}.
\end{equation}
Thus the zeroes of $H$ are exactly the points of $\Sigma $ where
$C=\pm 1/2$.

Next we explode the  equation (\ref{ya_esfera}) derivating it with
respect to $z$. To do that, we use again $\langle \xi, P \delta
\Phi \rangle =2i\,C \langle \xi, J \delta \Phi \rangle$, for any
normal vector field $\xi$, and the properties of $P$. This fact
gives that
\[
\partial_z \langle \delta \Phi, P \delta \Phi \rangle=2\langle\delta\delta\Phi,P\delta\Phi\rangle=
4(u_z\langle \delta \Phi, P \delta \Phi \rangle + i C \langle
\delta \delta \Phi, J \delta \Phi \rangle ),
\]
that joint with (\ref{ya_esfera}) produces
\begin{equation}\label{G_z}
\langle H,J\delta\Phi\rangle \, \partial_{z}\langle
H,J\delta\Phi\rangle=\frac{c_1-c_2}{4} \left(
u_z\langle\delta\Phi,P\delta\Phi\rangle+iC\langle\delta\delta\Phi,J\delta\Phi\rangle
\right).
\end{equation}

Now we define a new Hopf quadratic differential outside the zeroes of $H$
(that is, outside the zeroes of $\langle H, J\delta \Phi \rangle$) given by
\begin{equation}\label{Xi}
\Xi (z) =  \left( \frac{\langle \delta \delta \Phi, J \delta \Phi
\rangle}{\langle H, J\delta \Phi \rangle} +
\frac{c_1+c_2}{c_1-c_2} \langle \delta\Phi, P\delta\Phi
\rangle\right) (dz)^2.
\end{equation}
In order to prove that $\Xi$ is meromorphic, we use the Codazzi
equation of the immersion $\Phi$ and (\ref{C}) to get
\begin{equation}\label{f_zbarra}
\begin{array}{c}
\, \partial_{\bar z}\langle \delta \delta \Phi, J \delta \Phi
\rangle = \frac{1}{2} e^{2u} \partial_z \langle H, J\delta \Phi
\rangle \\ \\-u_z e^{2u} \langle H, J\delta \Phi \rangle -
\frac{i}{4} (c_1+c_2)e^{2u}C \langle \delta \Phi, P \delta \Phi
\rangle ,
\end{array}
\end{equation}
where we have used that $4\,\bar{R} (\bar{\delta}\Phi, \delta
\Phi,\delta \Phi,$ $ J\delta \Phi)=-i(c_1+c_2)e^{2u}C \langle
\delta \Phi, P \delta \Phi \rangle$. Using (\ref{G_z}) and
(\ref{f_zbarra}) we arrive at
\begin{equation}\label{+}
\begin{array}{c}
8\langle H, J\delta \Phi \rangle \partial_{\bar z} \,\langle
\delta\delta \Phi, J \delta \Phi \rangle =\\  \\
i\,e^{2u} C \left( (c_1-c_2) \langle \delta\delta \Phi, J \delta
\Phi \rangle - 16 \frac{c_1+c_2}{c_1-c_2} \langle H, J\delta \Phi
\rangle ^3 \right).
\end{array}
\end{equation}
It is easy to check now,  using (\ref{G_zbarra}) and (\ref{+}),
that the quadratic differential $\Xi$ is meromorphic. But if $p$
is a zero of $\langle H, J\delta \Phi \rangle$, (\ref{+}) says
that $p$ is also a zero of $\langle \delta\delta \Phi, J \delta
\Phi \rangle$ since at this point $C^2(p)=1/4$ from (\ref{|H|oC}).
In this way, $\Xi $ has no poles and so it is a holomorphic
differential. Hence, as $\Sigma$ is a sphere, necessarily $\Xi
\equiv 0$ and using (\ref{ya_esfera}) we obtain that
\begin{equation}\label{esfera2}
\langle\delta\delta\Phi,J\delta\Phi\rangle=-8\frac{c_1+c_2}{(c_1-c_2)^2}\langle
H,J\delta\Phi\rangle^3.
\end{equation}

Finally, taking modules in (\ref{esfera2}) and using (\ref{EqG})
and (\ref{|G|y|h|}) we obtain the following formula for the Gauss
curvature of our surface:
\begin{equation}\label{Ko|H|}
K=(c_1+c_2)C^2+\frac{|H|^2}{2}-8\frac{(c_1+c_2)^2}{(c_1-c_2)^4}|H|^6.
\end{equation}

\subsection{Geometric meaning of the vanishing of the Hopf differentials}
We are going to study the non trivial function $|H|$. To do that,
we define the complex function $h=2e^{-u}\langle
H,J\bar\delta\Phi\rangle$, which satisfies $|H|=|h|$. But using
(\ref{ya_esfera}) and (\ref{esfera2}),  we obtain that
\[
 h_{\bar{z}}=\left(u_{\bar z}+i\frac{c_1+c_2}{c_1-c_2}Ce^{u}h\right)h,
\]
which means that $h$ is locally the product of a positive function
and a holomorphic function. In fact, if $t$ is a
local solution of the equation
\[
t_{\bar{z}}=u_{\bar{z}}+i\frac{c_1+c_2}{c_1-c_2}Ce^{u}h,
\]
then $h=e^{t}(e^{-t}h)$, where $e^{t}$ is positive and from the
above equations $e^{-t}h$ is a holomorphic function.

Since $\Phi $ can not be minimal (see Proposition \ref{Prop2}),
using Lemma 4.1 in \cite{EGT}, we have that
\begin{equation}\label{EGT}
\int_\Sigma \triangle\log |H| dA = -2\pi N(|H|),
\end{equation}
 where $N(|H|)$ is
the sum of all orders for all zeroes of $|H|$.
From (\ref{|H|oC}), we get that
\begin{equation}\label{lalog1}
\triangle\log |H| = -\frac{2}{1-4C^2} \left( C\triangle C
+\frac{1+4C^2}{1-4C^2} |\nabla C|^2 \right).
\end{equation}
Hence, our next purpose is the computation of $\nabla C$ and
$\Delta C$.
To do that,
using (\ref{Frenet3}) and (\ref{C}) again, we deduce
\begin{equation}\label{C_z}
i e^{2u} C_z =2e^{-2u}\langle \bar{\delta} \Phi, P\bar{\delta }
\Phi \rangle \langle \delta \delta \Phi, J \delta \Phi \rangle
-\langle \delta \Phi, P\delta \Phi \rangle \, \langle H,
J\bar{\delta } \Phi \rangle.
\end{equation}
When we put in (\ref{C_z}) the information of (\ref{ya_esfera})
and (\ref{esfera2}) we arrive at
\[  C_z = \frac{2i}{c_1-c_2}|H|^2  \langle
H, J\delta \Phi \rangle \left( 1+4\frac{c_1+c_2}{(c_1-c_2)^2}
|H|^2 \right). \]
  Finally, using that $|\nabla C|^2=4e^{-2u}
|C_z|^2$ and (\ref{|G|y|h|}), we conclude that the modulus of the
gradient of $C$ is given by
\begin{equation}\label{gradC}
|\nabla C|^2=\frac{(1-4C^2)|H|^2}{4}\left(
1+4\frac{c_1+c_2}{(c_1-c_2)^2} |H|^2 \right)^2.
\end{equation}
As $\Delta C = 4e^{-2u}C_{z\bar z}$, we derivate (\ref{C_z}) with
respect to $\bar z$. After a long straightforward computation,
using (\ref{h_z}) and (\ref{|G|y|h|}) joint to (\ref{ya_esfera})
and (\ref{esfera2}), we deduce the formula for the Laplacian of
$C$, which is given by
\begin{equation}\label{deltaC}
\triangle C= -2|H|^2 C \left( 1+4\frac{c_1+c_2}{(c_1-c_2)^2} |H|^2
\right)^2 .
\end{equation}
Now we put (\ref{gradC}) and (\ref{deltaC}) in (\ref{lalog1}) to
reach
\begin{equation}\label{lalog2}
\triangle\log |H| = -\frac{|H|^2}{2}\left(
1+4\frac{c_1+c_2}{(c_1-c_2)^2} |H|^2 \right)^2 .
\end{equation}
Finally, using (\ref{|H|oC}) and (\ref{Ko|H|}) in (\ref{lalog2})
we get that
\[
\triangle\log |H|=K-|H|^2-\frac{c_1+c_2}{4}.
\]
Integrating this equation, using (\ref{EGT}) and the Gauss-Bonnet
and Poincaré-Hopf theorems, we conclude that
\begin{equation}\label{equalityin}
\int_{\Sigma}\left(|H|^2+\frac{c_1+c_2}{4}\right)dA=8\pi.
\end{equation}

This important geometric property of our Hamiltonian stationary
Lagrangian sphere is the key to prove the main result. In fact, we
finish this section by proving that $c_1$ and $c_2$ must be
positive and hence
$\Sigma_1\times\Sigma_2=\s^2_{c_1}\times\s^2_{c_2}$.

From (\ref{|H|oC}), we get that $|H|^2\leq\frac{c_1-c_2}{4}$ and
then (\ref{equalityin}) says that
\[
8\pi\leq \frac{c_1}{2}\hbox{Area}\,(\Sigma),
\]
which implies that $c_1>0$.

We now prove that not only $c_1>0$ but also $c_2>0$. The reasoning
starts from the Bochner formula for the tangent vector field $JH$:
\[
0=\int_\Sigma (K|JH|^2+|\nabla JH|^2-(\mathrm{div}\,JH)^2)dA.
\]
Using that $\mathrm{div}JH=0$ and the Lagrangian character of our
surface, the above formula becomes in
\begin{equation}\label{Bochner}
0=\int_\Sigma (K|H|^2+|\nabla^{\perp} H|^2)dA.
\end{equation}
Once more, the information about our immersion given in the equations
(\ref{h_z}), (\ref{ya_esfera}), (\ref{|H|oC}) and (\ref{+})
 allow us to express  $|\nabla^{\perp}H|^2$ in terms of $|H|^2$ and the associated Jacobian $C$.
After a non difficult long computation we arrive at
\[
|\nabla^{\perp}H|^2=2C^2 \left( \frac{(c_1+c_2)^2}{(c_1-c_2)^2}
|H|^4-\frac{1}{16}(c_1 -c_2)^2 \right).
\]
This formula joint with (\ref{Ko|H|}) becomes the integrand of
(\ref{Bochner}) in the following polynomial in $|H|^2$:
\begin{equation}\label{poly}
K|H|^2+|\nabla^\perp H|^2 = a |H|^4 + b|H|^2 + c,
\end{equation}
where
\[
a=\frac{4(c_1+c_2)^2 C^2-2c_1c_2}{(c_1-c_2)^2}, \  b=(c_1+c_2)C^2,
\ c =\frac{(c_1-c_2)^2}{8}.
\]
Suppose now that $c_2 \leq 0$. It is clear this implies
$a\geq 0$ and that the discriminant of the second degree polynomial
of (\ref{poly}) satisfies
$$b^2-4ac=\left(c_1 c_2- (c_1+c_2)^2 C^2\right)C^2 \leq 0.$$ Hence we
can deduce from (\ref{poly}) that $K|H|^2+|\nabla^\perp H|^2 \geq
0$. But then (\ref{Bochner}) gives that $K|H|^2+|\nabla^\perp H|^2
\equiv 0$ and therefore $$b^2-4ac=(c_1 c_2- (c_1+c_2)^2 C^2)C^2
\equiv 0$$ too. We have arrived at $C\equiv 0$ and this is
impossible according to the proof of Proposition \ref{Prop2}. In
this way we have proved that $c_2 $ must also be positive.

As a final conclusion of this section, {\it the Hamiltonian
stationary Lagrangian sphere $\Sigma $ lies in
$\Sigma_1\times\Sigma_2=\s^2_{c_1}\times\s^2_{c_2}$ and satisfies
\[
\int_{\Sigma}\left(|H|^2+\frac{c_1+c_2}{4}\right)dA=8\pi.
\]}

\subsection{Lagrangian spheres in $\s^2_{c_1}\times\s^2_{c_2}$ and the Willmore functional}

In this section  we are going to prove Theorem \ref{Th2} and
conclude the proof of the main Theorem \ref{Th1}.

We
consider $\s^2_{c_1}\times\s^2_{c_2}\subset\r^6$ and let
$\tilde{H}$ be the mean curvature vector of
$\Phi=(\phi,\psi):\Sigma \rightarrow \r^6$. Using (\ref{1=bo}) we
get that
\begin{equation}\label{Htilde}
\tilde{H}=H-\frac{1}{2}(c_1 \phi,c_2 \psi)
\end{equation}
 and hence
$|\tilde{H}|^2=|H|^2+\frac{c_1+c_2}{4}$. Thus the Willmore
functional $W(\Phi)$ of $\Sigma $ in $\r^6$ is given by
\begin{equation}\label{Willmore}
W(\Phi):=\int_{\Sigma}|\tilde{H}|^2\, dA=
\int_{\Sigma}\left(|H|^2+\frac{c_1+c_2}{4}\right)dA.
\end{equation}

We also need a result of Simon  and Li-Yau, that we recall now.
\vspace{0.1cm}
\begin{theorem}({\cite{S},\cite{LY}}\label{ThA})
Let
$\Phi:\Sigma\rightarrow\r^n$ be an immersion of a compact surface
$\Sigma$ with mean curvature vector $\tilde{H}$ and maximum
multiplicity $\mu$. Then
\[
\int_{\Sigma}|\tilde{H}|^2dA\geq 4\pi\mu,
\]
and the equality holds if and only if  $\tilde{H}$ is given on
$\tilde{\Sigma}=\Sigma-\{p_1,\dots,p_{\mu}\}$ by
$\tilde{H}=\frac{-2(\Phi-a)^{\perp}}{|\Phi-a|^2}$, where $\perp$
stands for normal component and $a=\Phi(p_i)$, for
all $1\leq i\leq \mu$. This condition about the mean
curvature $\tilde{H}$ means that
$\frac{\Phi-a}{|\Phi-a|^2}:\tilde{\Sigma}\rightarrow\r^n$ is a
complete minimal immersion with $\mu$ planar and embedded ends and
finite total curvature.
\end{theorem}

From Corollary 1, our Lagrangian sphere $\Sigma $ cannot be embedded
and so we have
that $\mu \geq 2$. Thus applying Theorem \ref{ThA} and using
(\ref{Willmore}) we get the inequality of Theorem \ref{Th2}.

Next we analyze the case of equality which, in particular, happens
to the Hamiltonian stationary Lagrangian sphere $\Sigma $ of
Theorem 1 according to (\ref{equalityin}).

As the equality holds, using Theorem \ref{ThA} we have that $\mu =
2$, $\Phi(p_1)=\Phi(p_2)=a=(a_1,a_2)\in \s^5_c$,
$c=\frac{c_1+c_2}{c_1 c_2} $, and the mean curvature vector
$\tilde{H}$ is given by $
\tilde{H}=\frac{-2(\Phi-a)^{\perp}}{|\Phi-a|^2}$, that joint with
(\ref{Htilde}) says
\[
a^{\perp}=\frac{|\Phi-a|^2}{2}(H-\frac{1}{2}(c_1\phi,c_2\psi))+\Phi.
\]
This property allows to prove that $\Phi$ lies in a affine
hyperplane of $\r^6$. In fact, if $\hat{a}=(-a_1,a_2)$, we have
that
\begin{equation}\label{last}
\langle \Phi, \hat{a} \rangle = \langle (-\phi,\psi),a \rangle =
\langle (-\phi,\psi),a^\perp
\rangle=\langle(-\phi,\psi),\Phi\rangle=1/c_2-1/c_1
\end{equation}
because $(-\phi,\psi)$ is a normal vector field to
$\s^2_{c_1}\times\s^2_{c_2}$ in $\r^6$.

Using again Theorem \ref{ThA}, we obtain a complete minimal
immersion
\[
\widehat{\Phi}=\frac{\Phi-a}{|\Phi -a|^2}:\Sigma \setminus \{
p_1,p_2 \} \rightarrow \r^6 ,
\]
with two planar ends and total curvature $-4\pi$. It is clear that
$\Phi$ can be recuperated from $\hat{\Phi}$ by
\begin{equation}\label{hatPhi1}
\Phi =a+\widehat{\Phi}/|\widehat{\Phi}|^2.
\end{equation}
Let us analyze in depth the immersion $\widehat{\Phi}$. First,
equation (\ref{last}) satisfying $\Phi$ becomes in
\begin{equation}\label{hatPhi2}
 \langle \widehat{\Phi}, \hat{a} \rangle= 0.
\end{equation}
On the other hand, using that $|\Phi -a|^2 =2(1/c-\langle \Phi, a
\rangle)$, it is clear that
\begin{equation}\label{caeR5}
\langle \widehat{\Phi}, a \rangle= -\frac{1}{2}.
\end{equation}
Hence from (\ref{hatPhi2}) and (\ref{caeR5}) we have deduced that
$\widehat{\Phi}$ lies in a 4-dimensional affine subspace of
$\r^6$. But then $\Hat{\Phi}$ is one of the embedded complex
surfaces of $\c^2$ with finite total curvature $-4\pi$ given in
Proposition 6.6 in \cite{HO} which, up to congruences and
dilations, can be described by the one-parameter family of
Lawlor's minimal cylinders (\cite{L}),
\[
\mathcal{C}_{t}=\{(z,w)\in\c^2\,/\,|z|^2-|w|^2=-1,\, \Re(zw)=\sinh
t \cosh t \},\ t \in \r.
\]
We now make use of the parametrization $(F_{t},G_{t}):\r
\times \s^1 \rightarrow \c^2$ of the cylinders
$\mathcal{C}_{t}$ obtained in \cite{CCh}, which are  given by
\begin{equation}\label{FG}
\begin{array}{c}
F_{t}(s_1,s_2)=\sqrt{c_1-c_2}  \, (\sinh t \cosh s_1 + i \,
\cosh t \sinh s_1)\, e^{i\,s_2} ,  \\  \\
G_{t}(s_1,s_2)=\sqrt{c_1-c_2} \, (\cosh t \cosh s_1 + i \,
\sinh t \sinh s_1)\, e^{-i\,s_2}.
\end{array}
\end{equation}
We can choose without restriction
$a_i=(0,1/\sqrt{c_i})\in\s^2_{c_i}$, $i=1,2$. From (\ref{hatPhi2})
and (\ref{caeR5}) we have that $\widehat{\Phi}$ is congruent to
some $\widehat{\Phi}_{t}$ of the one-parameter family of minimal
embeddings
\begin{equation}\label{invertida}
\widehat{\Phi}_{t}=\left(
F_{t},-\frac{\sqrt{c_1}}{4},G_{t},-\frac{\sqrt{c_2}}{4} \right):\r
\times\s^1\rightarrow\r^6\equiv \r^3 \times \r^3,\ t\in\r.
\end{equation}
Using then (\ref{invertida}), (\ref{FG}) and (\ref{hatPhi1}) joint
to the conformal map
\[
\begin{array}{c}
\r \times \s^1 \rightarrow \s^2 \\  \\
(s_1,e^{i\,s_2}) \mapsto (z,x)=\left(\frac{e^{i\, s_2}}{\cosh
s_1},\tanh s_1 \right),
\end{array}
\]
after a long straightforward computation, we conclude that our
original Lagrangian immersion $\Phi$ is congruent to some $\Phi_t
$  of the one-parameter family of Lagrangian immersions, $\Phi_t
=(\phi_t,\psi_t): \s^2 \rightarrow \s^2_{c_1} \times \s^2_{c_2}$,
$t \in \r$, given by
\begin{equation}\label{Phis}
\begin{array}{c}
\phi_t (z,x)=\frac{\textstyle 2\sqrt{c_1-c_2}}{\textstyle
(c_1 c_t^2 -c_2 s_t^2) + (c_1 s_t^2 -c_2
c_t^2)\, x^2} \times \\ \\
\left( ( s_t  + i \,  c_t \,x) z, \frac{\textstyle
(c_1c_t^2 -2c_1 -c_2 s_t^2)+(c_1 s_t^2
+2c_1-c_2 c_t^2)\, x^2}{\textstyle 2\sqrt{c_1}\sqrt{c_1-c_2}} \right), \\ \\
\psi_t (z,x)=\frac{\textstyle 2\sqrt{c_1-c_2}}{\textstyle
(c_1 c_t^2 -c_2 s_t^2) + (c_1 s_t^2 -c_2
c_t^2) \, x^2} \times \\  \\
\left( ( c_t  + i \,  s_t \,x) \bar{z}, \frac{\textstyle
(c_1c_t^2 -2c_2 -c_2 s_t^2)+(c_1 s_t^2 +2c_2-c_2
c_t^2)\, x^2}{\textstyle 2\sqrt{c_2}\sqrt{c_1-c_2}} \right),
\end{array}
\end{equation}
where $s_t=\sinh t$ and $c_t= \cosh t$.

It is now an exercise to check that $\Phi_{t},\, t\in\r$, satisfy
the equality (\ref{equalityin}), are embeddings except at the
poles of $\s^2$ where they have a double point and $\Phi_{-t}$ is
congruent to $\Phi_{t}$ for each $t>0$. This finishes the proof of
Theorem \ref{Th2}.

The proof of Theorem \ref{Th1} concludes when we compute the
divergence of the tangent vector field $JH_t$ for each immersion
$\Phi_{t}$, obtaining
\[
(\mathrm{div}\, JH_t )(z,x)= \frac{(c_2-c_1)\,(\sinh 2t)\,
x}{2(1+x^2)},\quad (z,x)\in\s^2,
\]
that proves that the immersion $\Phi_0$ given in (\ref{Phi_0}) is
the only Hamiltonian stationary one in the family.

\section{Stability properties of the sphere $\mathcal{S}_0$}
In this section we are interested in stability properties of the
Hamiltonian stationary Lagrangian immersion $\Phi_0$. We first
note that the Lagrangian variation $\Phi_{t},\, t\in\r$, of
$\Phi_0$ is really a Hamiltonian variation. In fact, it is not
difficult to check that the normal component of the variation
vector field $\frac{d\Phi_t}{dt}\left|_{_{t=0}} \right.
 $ is $J\nabla f$, where $f:\s^2\rightarrow\r$ is the
function given by
\[
f(z,x)=\frac{2(c_1 -c_2)}{c_1 c_2} \left[ \frac{(c_1 - c_2) x}{c_1
- c_2 x^2}
  - \frac{c_1 + c_2}{\sqrt{c_1 c_2}} \arctanh\left( \frac{\sqrt{c_2}}{\sqrt{c_1}}
  x\right)\right].
\]

On the other hand, a longer but easy computation says that the
area $A(t)$ of the induced metric on $\s^2$ by
$\Phi_{t}$ is given by
\[
    A(t) =
    \begin{cases}
        \frac{32 \pi}{s^2 - d^2} \left( s - \frac{2d^2}{\sqrt{s^2 - d^2}} \arctanh\left[ \frac{\sqrt{s - d}}{\sqrt{s + d}}\right]\right) & |t| < \frac{1}{2} \arccosh\left( \frac{c_1 + c_2}{c_1 - c_2} \right) \\
        \frac{64 \pi}{3(c_1 + c_2)} & |t |= \frac{1}{2} \arccosh\left( \frac{c_1 + c_2}{c_1 - c_2} \right) \\
        \frac{32 \pi}{s^2 - d^2} \left( s - \frac{2d^2}{\sqrt{d^2 - s^2}} \arctan\left[ \frac{\sqrt{d - s}}{\sqrt{d + s}}\right]\right) & |t | > \frac{1}{2} \arccosh\left( \frac{c_1 + c_2}{c_1 - c_2} \right) \\
    \end{cases}
    \]
where $s = c_1 + c_2$ and $d = (c_1 - c_2)\cosh (2t)$.

\begin{figure}[hb]\label{grafica}
\includegraphics[width=4cm]{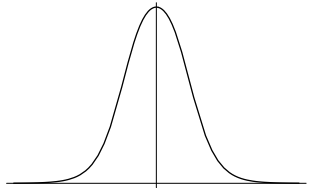}
\caption{Graphic of $A(t)$}
\end{figure}

 So it is easy to see that $A(t)$ has a maximum at
$t=0$ and, in addition, $\lim_{t\rightarrow
\infty}A(t)=0$ (see Figure 1). These two facts imply the following
conclusion:
\begin{quote}
{\it The Hamiltonian stationary Lagrangian immersion $\Phi_0$ is
unstable and there are no minimizers for the area in its
Hamiltonian isotopy class.}
\end{quote}

\end{document}